\documentclass[11pt, letterpaper, leqno]{amsart}
\usepackage{amsmath}
\usepackage{amsthm}
\usepackage{amsrefs}
\usepackage{lipsum}
\usepackage{amsopn}
\usepackage{amssymb}
\usepackage{enumerate} 
\numberwithin{equation}{section}
\usepackage{float}
\usepackage{graphicx}
\usepackage{titlesec}
\usepackage{fancyhdr}
\usepackage{fancyhdr}
\newtheorem*{condition1*}{Condition (1)}
\newtheorem*{condition2*}{Condition (2)}
\newtheorem*{condition3*}{Condition (3)}
\newtheorem{question}{Question}[section]

\newtheorem{lemma}{Lemma}[section]
\newtheorem{theorem}{Theorem}[section]

\theoremstyle{definition}

\DeclareMathOperator{\LOG}{Log}
\DeclareMathOperator{\Arg}{Arg}

\DeclareMathOperator{\arcsinh}{arcsinh}
\DeclareMathOperator{\dist}{dist}
\DeclareMathOperator{\arc}{arc}
\titleformat{\section}
  {\normalfont\scshape}{\thesection}{1em}{\centering}
\titleformat{\subsection}
{\normalfont\scshape}{\thesubsection}{1em}{\centering}

\begin{document}
\title{On a property of harmonic measure on simply connected domains}

\author{Christina Karafyllia}  
\address{Department of Mathematics, Aristotle University of Thessaloniki, 54124, Thessaloniki, Greece}
\email{karafyllc@math.auth.gr}   
\thanks{I would like to thank Professor D.\ Betsakos, my thesis advisor, for his advice during the preparation of this work and the Onassis Foundation for the scholarship I receive during my Ph.D.\ studies.\ I would also like to thank the referees for their useful remarks and their suggestions about simplifying some proofs.}

\fancyhf{}
\renewcommand{\headrulewidth}{0pt}

\fancyhead[RO,LE]{\small \thepage}
\fancyhead[CE]{\small On a property of harmonic measure on simply connected domains}% odd page header and number to right top
\fancyhead[CO]{\small Christina Karafyllia} 
\fancyfoot[L,R,C]{}


%\date{}

\begin{abstract} Let $D \subset \mathbb{C}$ be a domain with $0 \in D$.\ For $R>0$, let
${{\hat \omega }_D}\left( {R} \right)$ denote the harmonic measure of $ D \cap \left\{ {\left| z \right| = R} \right\}$ at $0$ with respect to the domain $
D \cap \left\{ {\left| z \right| < R} \right\}
$ and
${\omega _D}\left( {R} \right)$ denote the harmonic measure of $\partial D \cap \left\{ {\left| z \right| \ge R} \right\}$ at $0$ with respect to $D$.\ The behavior of the functions ${\omega _D}$ and ${{\hat \omega }_D}$ near $\infty$ determines (in some sense) how large $D$ is.\ However, it is not known whether the functions ${\omega _D}$ and ${{\hat \omega }_D}$ always have the same behavior when $R$ tends to $\infty$.\ Obviously, ${\omega _D}\left( {R} \right) \le {{\hat \omega }_D}\left( {R} \right)$ for every $R>0$.\ Thus, the arising question, first posed by Betsakos, is the following: Does there exist a positive constant $C$ such that for all simply connected domains $D$ with $0 \in D$ and all $R>0$,
	\[{\omega _D}\left( {R} \right) \ge C{{\hat \omega }_D}\left( {R} \right)?\]
In general, we prove that the answer is negative by means of two different counter-examples.\ However, under additional assumptions involving the geometry of $D$, we prove that the answer is positive.\ We also find the value of the optimal constant for starlike domains.
\end{abstract}

\maketitle
\section{Introduction}\label{section1}

We will give an answer to a question of Betsakos (\cite[p.\ 788]{Bet}) about a property of harmonic measure.\  For a domain $D$, a point $z \in D$ and a Borel subset $E$ of $\overline D $, let ${\omega _D}\left( {z,E} \right)$ denote the harmonic measure at $z$ of $\overline E$ with respect to the component of $D \backslash {\overline E}$ containing $z$.\ The function ${\omega _D}\left( { \cdot ,E} \right)$ is exactly the solution of the generalized Dirichlet problem with boundary data $\varphi  = {1_E}$   (see \cite[ch.\ 3]{Ahl}, \cite[ch.\ 1]{Gar} and \cite[ch.\ 4]{Ra}).\ The probabilistic interpretation of harmonic measure is that, given a domain $D$, a point $z\in D$ and a set $E \subset \partial D$, the harmonic measure ${\omega _D}\left( { z ,E} \right)$ is the probability that a Brownian motion started at $z$ will first hit the boundary of $D$ in the set $E$.

Let $D \subset \mathbb{C}$ be a domain with $0 \in D$.\ For $R>0$, we set
\[{\omega _D}\left( {R} \right) = {\omega _D}\left( {0,\partial D \cap \left\{ {z:\left| z \right| \ge R} \right\}} \right)\]
and
\[{{\hat \omega }_D}\left( {R} \right) = {\omega _D}\left( {0,
	 D  \cap 
	\left\{ {z:\left| z \right| = R} \right\}} \right).\]
The behavior of the functions ${\omega _D}$ and ${{\hat \omega }_D}$ near $\infty$ determines (in some sense) how large $D$ is and it has been studied from various viewpoints.\ For example, in \cite{Ts} and \cite[p.\ 111-118]{Tsu} Tsuji proved bounds for the growth of ${{\hat \omega }_D}\left( {R} \right)$ in terms of the size of the maximal arcs on $\left\{ {z:\left| z \right| = R} \right\}$.\ Tsuji's inequalities can be used to obtain estimates for the maximum modulus, means and coefficients of various classes of $p-$valent functions (see also \cite[ch.\ 8]{Hay}).\ In \cite{We} Hayman and Weitsman used ${{\hat \omega }_D}\left( {R} \right)$ to estimate the means and hence the coefficients of functions when information is known about their value distribution.\ With the aid of ${\omega _D}\left( {R} \right)$ and ${{\hat \omega }_D}\left( {R} \right)$, Sakai \cite{Sa} gave an integral representation of the least harmonic majorant of $\left| x \right| ^p$ in an open subset $D$ of $\mathbb{R}^n$ with $0 \in D$ and proved isoperimetric inequalities for it.\ Ess{\' e}n, Haliste, Lewis and Shea (\cite{Ess}, \cite{Esse}) also studied the problem of harmonic majoration in higher dimensions in terms of the geometry of $D$ by using ${\omega _D}\left( {R} \right)$ and ${{\hat \omega }_D}\left( {R} \right)$.\ In \cite[p.\ 1348]{Soly} Solynin proved an estimate of ${{ \omega }_D}\left( {R} \right)$ when $D = f\left( \mathbb{D} \right)$ and $f$ is in the class $S$ of functions which are regular and univalent in the unit disk and $f\left( 0 \right)=0$, $f'\left( 0 \right) = 1$.\ Baernstein \cite{Bae} proved an integral formula involving ${{\hat \omega }_D}\left( {R} \right)$ and Green's function.

In \cite{Es} Ess{\' e}n proved that every analytic function $f:\mathbb{D} \to D$ belongs to the Hardy space $H^p$ for some $p>0$ if and only if for some constants $q$ and $C$, we have ${{\hat \omega }_D}\left( R \right) \le C{R^{ - q}}$ for every $R \ge 1$.\ With the aid of Ess{\' e}n's result, Kim and Sugawa \cite{Kim} proved that the Hardy number, ${\rm {h}}\left( D \right)$, of a plane domain $D$ with $0 \in D$, can be determined by
\[{\rm {h}}\left( D \right) =  - \mathop {\lim \sup }\limits_{R \to  + \infty } \frac{{\log {{\hat \omega }_D}\left( R \right)}}{{\log R}}.\]

In \cite{Be} Betsakos studied another problem involving ${\omega _D}\left( {R} \right)$.\ Let $\mathcal{B}$ be the family of all simply connected domains $D \subset \mathbb{C}$ such that $0 \in D$ and there is no disk of radius larger than $1$ contained in $D$.\ It is obvious that if $D \in \mathcal{B}$ then ${\omega _D}\left( {R} \right)$ is a decreasing function of $R$.\ In fact, ${\omega _D}$ decays exponentially as it is proved that  there exist positive constants $\beta$ and $C$ such that ${\omega _D}\left( R \right) \le C{e^{ - \beta R}}$, for every $D \in \mathcal{B}$ and every $R>0$.\ The problem studied in \cite{Be} is to find the optimal exponent $\beta$. 
	 
Poggi-Corradini (see \cite[p.\ 33-34]{Co}, \cite{Co1}, \cite{Co2}) studied ${\omega _D}\left( {R} \right)$ and ${{\hat \omega }_D}\left( {R} \right)$ in relation with conformal mappings in Hardy spaces.\ In fact, if $D$ is an unbounded simply connected domain with $0 \in D$ and $\psi$ is a conformal mapping of $\mathbb{D}$ onto $D$, then he proved that
\[\psi \in {H^p}\left( \mathbb{D} \right) \Leftrightarrow \int_0^{ + \infty } {{R^{p - 1}}{\omega _D}\left( R \right)dR}  <  + \infty  \Leftrightarrow \int_0^{ + \infty } {{R^{p - 1}}{{\hat \omega }_D}\left( R \right)dR}  <  + \infty.\] 
To establish the last equivalence, Poggi-Corradini first proved that there exists a constant $M_0>1$ such that for all $R>0$,
\begin{equation}\label{co}
{\omega _D}\left( R \right) \ge \frac{1}{2}{{\hat \omega }_D}\left( {{M_0}R} \right).
\end{equation}

All the results mentioned above are some of the estimates and applications of ${\omega _D}$ and ${{\hat \omega }_D}$ that have been made over time.\ However, it is still unknown whether the functions ${\omega _D}$ and ${{\hat \omega }_D}$ always have the same behavior when $R$ tends to $\infty$.\ Obviously, by the maximum principle, for every $R>0$,
\[{\omega _D}\left( {R} \right) \le {{\hat \omega }_D}\left( {R} \right)\]
but all we know about the inverse inequality is (\ref{co}).\ Thus, a natural question, first posed in \cite[p.\ 788]{Bet} by Betsakos, is the following:
\begin{question}\label{con} Does there exist a positive constant $C$ such that for a class of domains $D$ (such as simply connected, starlike etc.) with $0 \in D$ and every $R>0$, 
\[{\omega _D}\left( {R} \right) \ge C{{\hat \omega }_D}\left( {R} \right)?\]
\end{question}
In this paper we prove that for simply connected domains the answer is negative by means of two different counter-examples.\ However, under additional assumptions involving the geometry of the domains, we prove that the answer is positive and we also find the value of the optimal constant for starlike domains.

\begin{figure}[H]
		\begin{center}
			\includegraphics[scale=0.55]{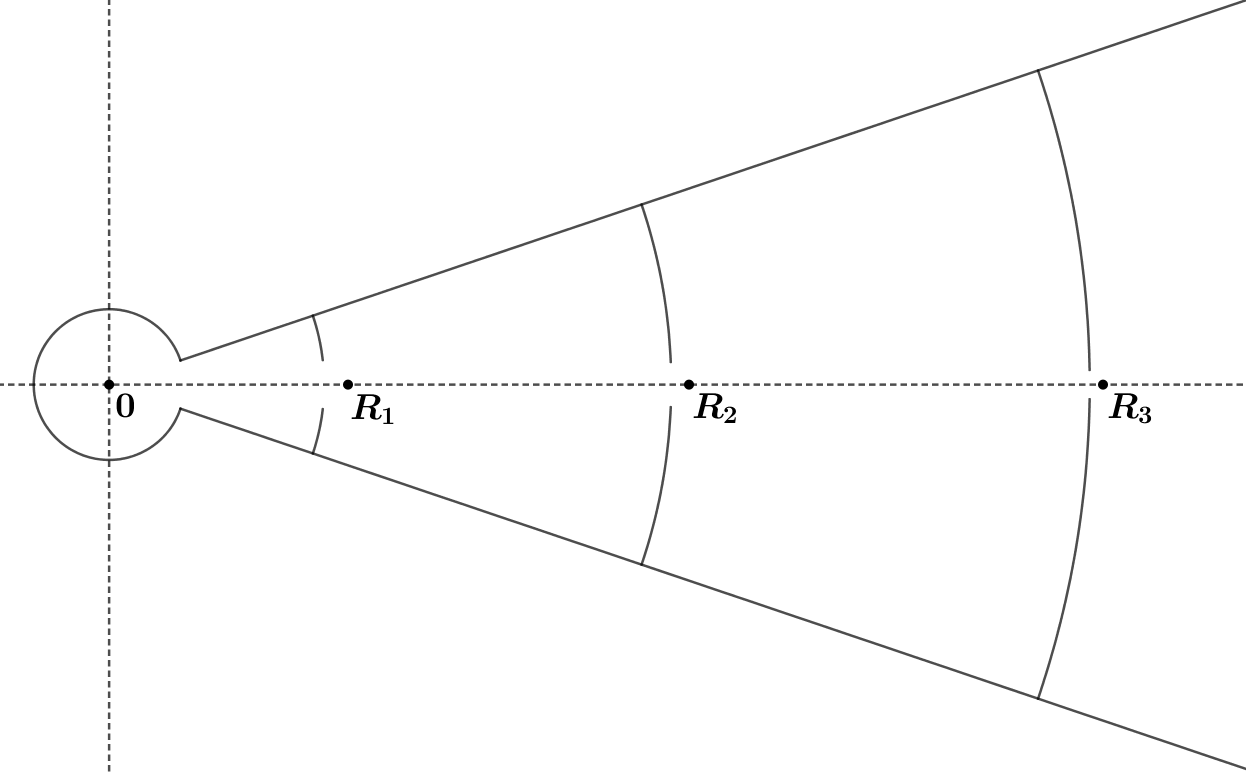}
			\vspace*{0cm}
			\caption{}
			\label{dom1}
		\end{center}
\end{figure}
 
In Section \ref{section3}, we construct the simply connected domain $D$ of Fig.\ \ref{dom1} and prove that there exists a sequence of positive numbers ${\left\{ {{R_n}} \right\}_{n \in \mathbb{N}}}$ such that 
\[\mathop {\lim }\limits_{n \to  + \infty } \frac{{{{\hat \omega }_D}\left( {{R_n}} \right)}}{{{\omega _D}\left( {{R_n}} \right)}} = + \infty,\]
which implies that there does not exist a positive constant $C$ such that
${\omega _D}\left( {R} \right) \ge C{{\hat \omega }_D}\left( {R} \right)$ for every $R>0$.\ As we see in the proof, this result is due to the fact that the hyperbolic distance between the point $R_n$ and the hyperbolic geodesic, $\Gamma _n$, joining the endpoints of the arc $ D  \cap \left\{ {\left| z \right| = R_n } \right\}$ in $D$ tends to infinity as $n\to +\infty$.\ In other words, there does not exist a positive constant $c$ such that $D  \cap \left\{ {\left| z \right| = R_n } \right\} \subset \left\{ {z \in D:d_D \left( {z,\Gamma _n } \right) < c} \right\}$ for every $n \in \mathbb{N}$.\ Note that $d_D \left( {z,\Gamma _n } \right)$ denotes the hyperbolic distance between $z$ and $\Gamma _n$ in $D$, which we define in Section \ref{section2}.\ Now we consider the following condition on the simply connected domain $D$:

\begin{condition1*}
There exists a constant $c>0$ such that, for every $R>0$, every arc of $D \cap \left\{ {z:\left| z \right| = R} \right\}$ lies in a hyperbolic $c$-neighborhood of the hyperbolic geodesic joining its endpoints. 	
\end{condition1*}

\begin{figure}[H] 
	\begin{center}
		\includegraphics[scale=0.55]{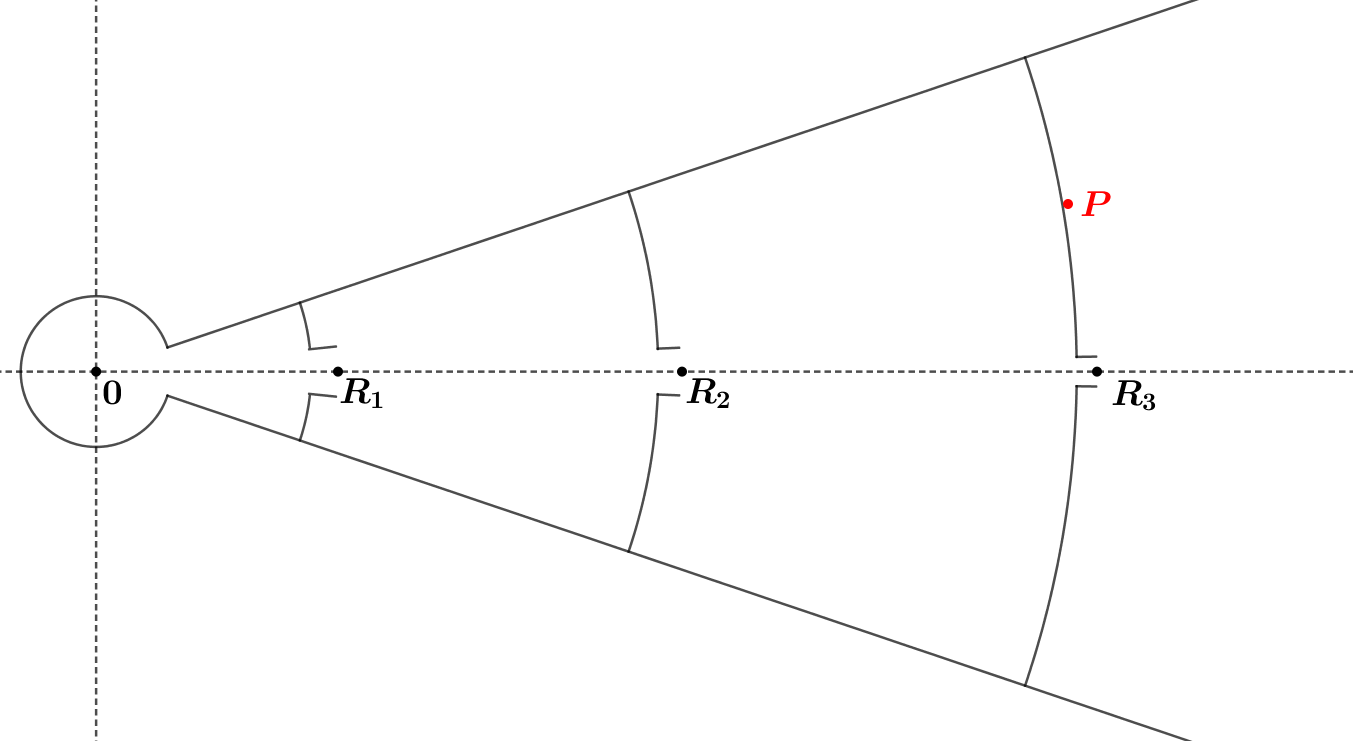}
		\caption{}
		\label{do2}
	\end{center}
\end{figure}

The arising question is whether the answer to the Question \ref{con} is positive for simply connected domains that satisfy Condition (1).\ However, we prove that this condition is not enough by constructing, in Section \ref{section4}, the simply connected domain $D$ of Fig.\ \ref{do2}, which comes from a small variation of the domain of Fig.\ \ref{dom1}.\ In fact, there exists a sequence of positive numbers ${\left\{ {{R_n}} \right\}_{n \in \mathbb{N}}}$ such that, despite the fact that Condition (1) is satisfied, we have again 
\[\mathop {\lim }\limits_{n \to  + \infty } \frac{{{{\hat \omega }_D}\left( {{R_n}} \right)}}{{{\omega _D}\left( {{R_n}} \right)}} = + \infty.\] 
This time, this is due to the fact that there exists a prime end $P$ of $\partial{D}$ that is inside the disk $\left\{ {z:\left| z \right| < R_n} \right\}$ but every arc in $D$ joining $0$ to $P$ intersects the circle $\left\{ {z:\left| z \right| = R_n} \right\}$.\ See, for example, the prime end $P$ in Fig.\ \ref{do2}.\ So, we consider the following condition:  
 
\begin{condition2*}
For every $R>0$, there does not exist any prime end $P$ of $\partial{D}$ that is inside the disk $\left\{ {z:\left| z \right| < R} \right\}$ but every arc in $D$ joining $0$ to $P$ intersects the circle $\left\{ {z:\left| z \right| = R} \right\}$.
\end{condition2*}

Note that in the first counter-example (Section \ref{section3}) Condition (2) is satisfied, since it is obvious that there do not exist such prime ends.\ These two counter-examples show that Conditions (1) and (2) are necessary if we want to give a positive answer to the Question \ref{con}.\ But are they enough? In Section \ref{section5}, we actually prove that if a simply connected domain satisfies Conditions (1) and (2), then there exists a positive constant $K=K\left( c \right)$ such that for every $R>0$,
\[{{\hat \omega }_D}\left( {R} \right) \le K{\omega _D}\left( {R} \right).\]
Moreover, we prove that we can find the value of this constant if we retain Condition (2) and replace Condition (1) with the following condition:
\begin{condition3*}
For every $R>0$ and for every arc of $D \cap \left\{ {z:\left| z \right| = R} \right\}$, the hyperbolic geodesic joining its endpoints lies entirely in $\overline{D} \cap \left\{ {z:\left| z \right| \le R} \right\}$.
\end{condition3*} 

So, having these results in mind, in Section \ref{section5}, we prove the theorem below which gives a positive answer to the Question \ref{con}. 

\begin{theorem}\label{kyrio} Let $D \subset \mathbb{C}$ be a simply connected domain with $0 \in D$.\ With the notation above, if Conditions $\rm{(1)}$ and $\rm{(2)}$ are satisfied, then there exists a positive constant $K=K\left( c \right)$ such that for every $R>0$,
	\[{{\hat \omega }_D}\left( {R} \right) \le K{\omega _D}\left( {R} \right).\]
If Conditions $\rm{(2)}$ and $\rm{(3)}$ are satisfied, then for every $R>0$,
\[{{\hat \omega }_D}\left( {R} \right) \le 2{\omega _D}\left( {R} \right).\]	
\end{theorem}

Finally, recall that a domain $D$ in $\mathbb{C}$ is called starlike with respect to $0$, if for every point $z \in D$, the segment of the straight line from $0$ to $z$, $\left[ {0,z} \right]$, lies entirely in $D$.\ In Section \ref{section6}, we prove that starlike domains satisfy Conditions (2) and (3) and that $2$ is the optimal constant: 

\begin{theorem}\label{star} Let $D$ be a starlike domain in $\mathbb{C}$.\ Then for every $R>0$,
	\[{{\hat \omega }_D}\left( {R} \right) \le 2{\omega _D}\left( {R} \right)\]
	and the constant $2$ is best possible.
\end{theorem}

In Section \ref{section2}, we introduce some preliminaries such as notions and results in hyperbolic geometry and basic properties of harmonic measure.\ In Sections \ref{section3} and \ref{section4}, we present the counter-examples of Fig.\ \ref{dom1} and \ref{do2} respectively, and in Sections \ref{section5} and \ref{section6}, we prove Theorems \ref{kyrio} and \ref{star} respectively. 

\section{Preliminary results}\label{section2}

\subsection{Results in hyperbolic geometry}

For the unit disk $\mathbb{D}$ the density of the hyperbolic metric is 
\[{\lambda _\mathbb{D}}\left( z \right) = \frac{2}{{1 - {{\left| z \right|}^2}}}.\]
Let $\Omega$ be a hyperbolic region in the complex plane $\mathbb{C}$; that is, $\mathbb{C}\backslash \Omega $ contains at least two points.\ If $f$ is a holomorphic universal covering projection of $\mathbb{D}$ onto $\Omega$ then the density $\lambda _{\Omega}$ is determined from
\[{\lambda _{\Omega}}\left( {f\left( z \right)} \right)\left| {f'\left( z \right)} \right| = \frac{2}{{1 - {{\left| z \right|}^2}}}\]
(see \cite[p.\ 236]{Mi}).\ The determination of $\lambda _{\Omega}$ is independent of the choice of the holomorphic covering projection onto $\Omega$.\ If $\Omega$ is simply connected, then $f$ is a conformal mapping of $\mathbb{D}$ onto $\Omega$.\ We note that in this paper we work on simlpy connected domains.

The hyperbolic distance between two points $z,w$ in $\mathbb{D}$ is defined by 
\[{d_\mathbb{D}}\left( {z,w} \right) = \log \frac{{1 + \left| {\frac{{z - w}}{{1 - z\bar w}}} \right|}}{{1 - \left| {\frac{{z - w}}{{1 - z\bar w}}} \right|}}\](see \cite[ch.\ 1]{Ahl}, \cite[p.\ 11-28]{Bea}).\ It is conformally invariant and thus it can be defined on any simply connected domain $D \ne \mathbb{C}$ as follows: If $f$ is a Riemann mapping of $\mathbb{D}$ onto $D$ and $z,w \in D$, then 
${d_D}\left( {z,w} \right) = {d_\mathbb{D}}\left( {{f^{ - 1}}\left( z \right),{f^{ - 1}}\left( w \right)} \right)$.
Also, for a set $E \subset D$, we define ${d_D}\left( {z,E} \right): = \inf \left\{ {{d_D}\left( {z,w} \right):w \in E} \right\}$.

The following theorem is known as Minda's reflection principle \cite[p.\ 241]{Mi}.\ First, we introduce some notation: If $\Gamma$ is a straight line (or circle), then $R$ is one of the half-planes (or the disk) determined by $\Gamma$ and ${\Omega ^*}$ is the reflection of a hyperbolic region $\Omega$ in $\Gamma$ .

\begin{theorem} \label{rp} Let $\Omega$ be a hyperbolic region in $\mathbb{C}$ and $\Gamma$ be a straight line or circle with $\Omega  \cap \Gamma  \ne \emptyset $.\ If $\Omega \backslash R \subset \Omega^*$, then 
	\[{\lambda _{{\Omega ^*}}}\left( z \right) \le {\lambda _\Omega }\left( z \right)\]
	for all $z \in \Omega \backslash \overline R $.\ Equality holds if and only if $\Omega$ is symmetric about $\Gamma$.
\end{theorem}
\noindent
A generalization of Theorem \ref{rp} was proved by Solynin in \cite{Sol}. 

\subsection{Quasi-hyperbolic distance}

The hyperbolic distance between $z_1,z_2 \in D$ can be estimated by the quasi-hyperbolic distance, ${\delta _D}\left( {z_1,z_2} \right)$, which is defined by
\[{\delta _D}\left( {{z_1},{z_2}} \right) = \mathop {\inf }\limits_{\gamma :{z_1} \to {z_2}} \int_\gamma  {\frac{{\left| {dz} \right|}}{{d\left( {z,\partial D} \right)}}}, \]
where the infimum ranges over all the paths  connecting $z_1$ to $z_2$ in $D$ and $d\left( {z,\partial D} \right)$ denotes the Euclidean distance of $z$ from $\partial D$.\ Then it is proved that $\left( {{1 \mathord{\left/
			{\vphantom {1 2}} \right.
			\kern-\nulldelimiterspace} 2}} \right){\delta _D} \le {d_D} \le 2{\delta _D}$ (see \cite[p.\ 33-36]{Bea}, \cite[p.\ 8]{Co}).

\subsection{Harmonic measure}

If $E \subset {{\overline {\mathbb{D}}\backslash \left\{ 0 \right\}}}$, then a special case of the Beurling-Nevanlinna projection theorem (see \cite[p.\ 43-44]{Ahl}, \cite[p.\ 105]{Gar} and \cite[p.\ 120]{Ra}) is the following: 

\begin{theorem}\label{bene} Let $E \subset {{\overline {\mathbb{D}}\backslash \left\{ 0 \right\}}}$ be a closed, connected set intersecting the unit circle.\ If ${r_0} = \min \left\{ {\left| z \right|:z \in E} \right\}$ and ${E^ * } = \left\{ { - \left| z \right|:z \in E} \right\} = \left( { - 1,} \right.\left. { - {r_0}} \right]$, then
	\[{\omega _{\mathbb{D}}}\left( {0,E} \right) \ge {\omega _{\mathbb{D} }}\left( {0,{E^ * }} \right) = \frac{2}{\pi }\arcsin \frac{{\left( {1 - {r_0}} \right)}}{{\left( {1 + {r_0}} \right)}}.\]
\end{theorem}

Next theorem states the strong Markov property for harmonic measure, which follows from the probabilistic interpretation of harmonic measure (see \cite[p.\ 282]{Be} and \cite[p.\ 88]{Port}).

\begin{theorem}\label{markov}
Let $D_1$ and $D_2$ be two domains in $\mathbb{C}$.\ Assume that $D_1 \subset {D_2} $ and let	$F \subset \partial{D_2}$ be a closed set.\ If $\sigma=\partial{D_1}\backslash \partial {D_2}$, then for $z \in {D_1}$,
\[{\omega _{{D_2}}}\left( {z,F} \right) = {\omega _{{D_1}}}\left( {z,F} \right) + \int_\sigma  {{\omega _{{D_1}}}\left( {z,ds} \right){\omega _{{D_2}}}\left( {s,F} \right)}. \]
\end{theorem} 

The following result of Balogh and Bonk \cite{Bo} gives an estimate of the logarithmic capacity of a set $E\subset \partial \mathbb{D}$.\ But this also proves an estimate of harmonic measure because if $E$ is a finite union of closed arcs in $\partial \mathbb{D}$, then $\omega _\mathbb{D} \left( {0,E} \right) \le {\rm cap}E$ (see \cite[p.\ 164]{Gar}).

\begin{theorem}\label{bonk}
There exists a universal constant $K>0$ with the following property.\ Suppose $f:\mathbb{D} \to \mathbb{C}$ is a conformal mapping with $\dist\left( {f\left( 0 \right),\partial f\left( \mathbb{D} \right)} \right)\\= d$.\ If $E_f \left( R \right)$ is the set of all $\zeta \in \partial \mathbb{D}$ with $\rm{length}$ $f\left( {\left[ {0,\zeta } \right)} \right) \ge R > 0$, then
\[{\rm cap}E_f \left( R \right) \le K\sqrt {\frac{d}{R}}. \]

\end{theorem}

Next theorem states a relation between harmonic measure and hyperbolic distance, which we prove in \cite{Kara}.

\begin{theorem}\label{geod}
	Let $\Gamma $ be the hyperbolic geodesic joining two points $z_1,z_2 \in \partial{\mathbb{D}}$ in $\mathbb{D}$.\ Then
	\[{e^{ - {d_{\mathbb{D}}}\left( {0,\Gamma } \right)}} \le {\omega _{\mathbb{D}}}\left( {0,\Gamma } \right) \le \frac{4}{\pi }{e^{ - {d_{\mathbb{D}}}\left( {0,\Gamma } \right)}}.\]
\end{theorem}

\section{First counter-example}\label{section3}

\begin{figure}[H] 
	\begin{center}
		\includegraphics[scale=0.55]{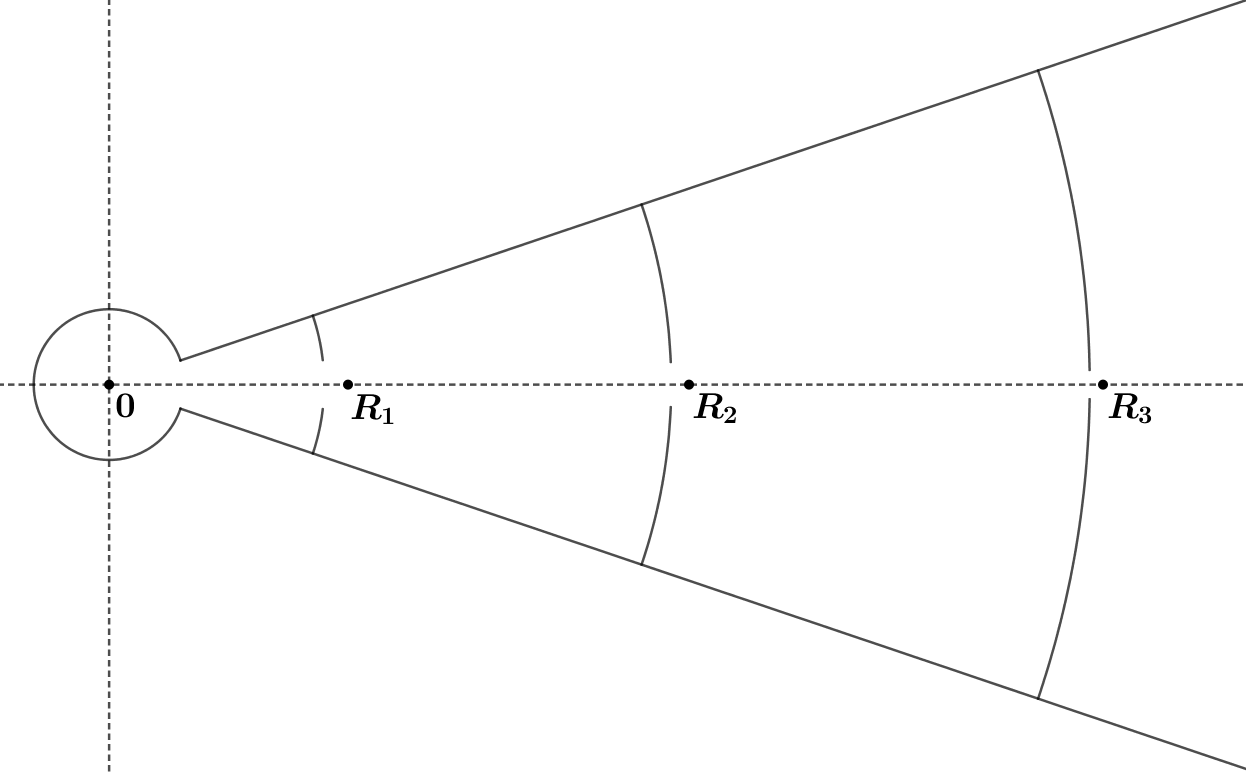}
		\caption{The simply connected domain $D$.}
		\label{dom}
	\end{center}
\end{figure}

Hereinafter, we use the notation $D\left( {z,r} \right) := \left\{ {w \in \mathbb{C}:\left| {w - z} \right| < r} \right\}$ for some $z \in \mathbb{C}$ and some $r>0$.\ Let $D$ be the simply connected domain of Fig.\ \ref{dom}, namely,
\[D = \mathbb{D} \cup \left( \left\{ {z \in \mathbb{C} :\left| {\Arg{z}} \right| < 1} \right\}\backslash \bigcup\limits_{n = 1}^{ + \infty } {\left\{ {z \in \partial D\left( {0,{e^n}} \right):\frac{1}{{{40^n}}} \le \left| {\Arg{z}} \right| \le 1} \right\}}\right)\]
and consider the sequence ${\left\{ {{R_n}} \right\}_{n \in \mathbb{N}}}$ with ${R _n} = {e^{n + \frac{1}{{{{40}^n}}}}}$ for every $n \in \mathbb{N}$.

\begin{theorem}\label{ant1} With the notation above, the simply connected domain $D$ has the following properties:
\begin{enumerate}[\rm(i)]
\item $D$ satisfies Condition {\rm (2)}.
\item $D$ does not satisfy Condition {\rm (1)}. 
\item $\mathop {\lim }\limits_{n \to  + \infty } \frac{{{{\hat \omega }_D}\left( {{R_n}} \right)}}{{{\omega _D}\left( {{R_n}} \right)}} = + \infty.$	
\end{enumerate}	
\end{theorem}

\proof 
Property (i) is immediate by the construction of $D$.\ So, we prove properties (ii) and (iii) (for a similar calculation see \cite{Ka}).\ The Riemann mapping theorem implies that there exists a conformal mapping $\psi $ from $\mathbb{D}$ onto $D$ such that $\psi \left( 0 \right) = 0$.\ For $n \in \mathbb{N}$, we set ${F_{R_n}} = \left\{ {z \in \mathbb{D}:\left| {\psi \left( z \right)} \right| = {R_n}} \right\}$ and ${E_{R_n} } = \left\{ {\zeta  \in \partial \mathbb{D}:\left| {\psi \left( \zeta \right)} \right| \ge {R_n}} \right\}$.\ Also, for $n \in \mathbb{N}$, let $\Gamma_{R_n}$ be the hyperbolic geodesic joining the endpoints of $F _{{R _n}}$ in $\mathbb{D}$.

\begin{figure}[H] 
	\begin{center}
		\includegraphics[scale=0.6]{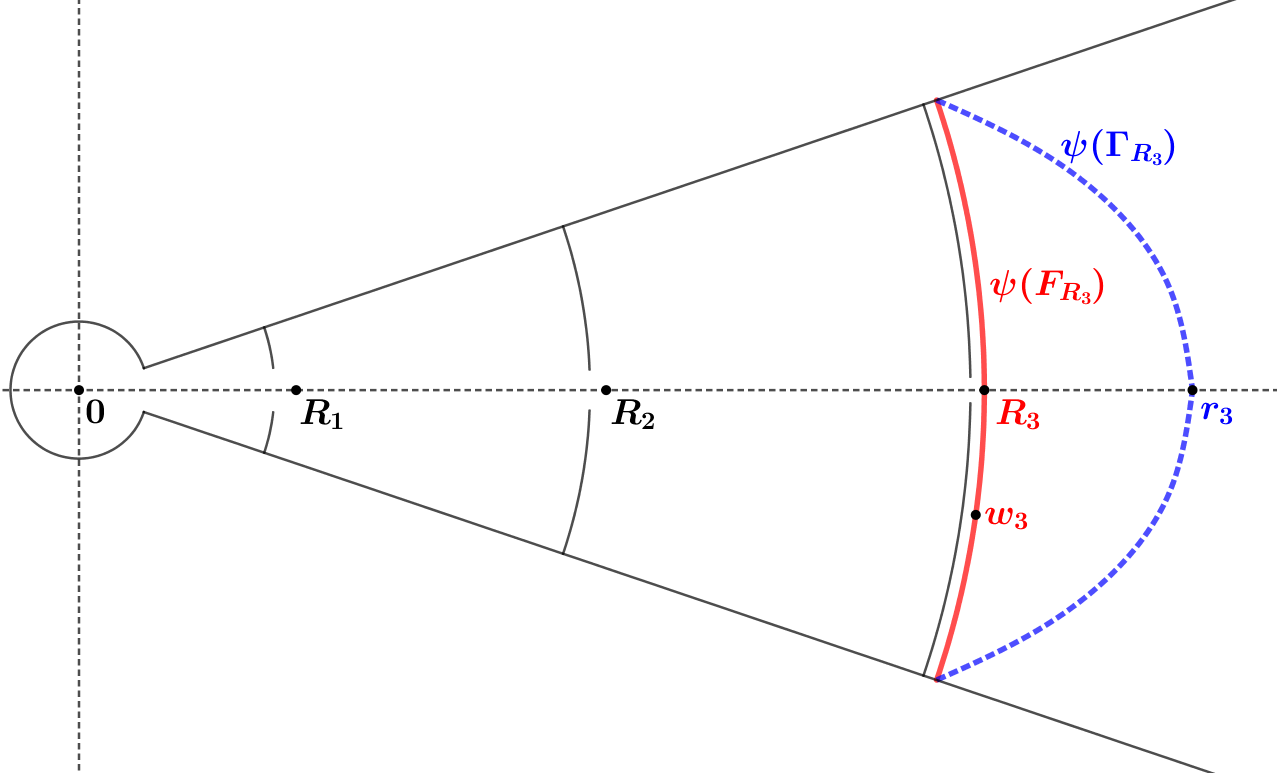}
		\caption{The crosscuts  $\psi \left( {{F_{{R _n}}}} \right)$ and $\psi \left( {{\Gamma_{{R _n}}}} \right)$ in case $n=3$.}
		\label{eik}
	\end{center}
\end{figure}

By Theorem \ref{bene} and the definition of hyperbolic distance we can easily infer that for every $n \in \mathbb{N}$,

\begin{equation}\nonumber
{\omega _\mathbb{D}}\left( {0,{F_{R_n}}} \right) \ge \frac{2}{\pi }{e^{ - {d_\mathbb{D}}\left( {0,{F_{R_n}}} \right)}}
\end{equation}
(see \cite[p.\ 35]{Co}).\ So, by the conformal invariance of harmonic measure and hyperbolic distance, we have
\begin{equation}\label{bn}
\hat \omega _D \left( {R_n } \right) = \omega _D \left( {0,\psi \left( {F_{R_n } } \right)} \right) \ge \frac{2}{\pi }e^{ - d_D \left( {0,\psi \left( {F_{R_n } } \right)} \right)}. 
\end{equation} 
Now fix a number $n>2$.\ If $z \in D$ and ${g_{D}}\left( {\cdot,\cdot
} \right)$ denotes the Green function for $D$ (see \cite[p.\ 41-43]{Gar}, \cite[p.\ 106-115]{Ra}), then
\[{d_{D}}\left( 0,z \right) = \log \frac{{1 + {e^{ - {g_{D}}\left( {0,z} \right)}}}}{{1 - {e^{ - {g_{D}}\left( {0,z} \right)}}}}\]
(see \cite[p.\ 12-13]{Bea} and \cite[p.\ 106]{Ra}).\ For every $w_n \in \psi \left( {{F_{{R _n}}}} \right)\backslash \left\{ {{R_n}} \right\}$ (see Fig.\ \ref{eik}), we infer, by a symmetrization result, that
\[{g_{D}}\left( {0,{R_n}} \right) \ge {g_{D}}\left( {0, w_n} \right)\]
(see Lemma 9.4 in \cite[p.\ 659]{Hay}).\ Since
\[f\left( x \right) = \log \frac{{1 + {e^{ - x}}}}{{1 - {e^{ - x}}}}\]
is a decreasing function on $\left( {0, + \infty } \right)$, we have that
\[{d_{D}}\left( {0, \psi \left( {{F_{{R _n}}}} \right)} \right)={d_{D}}\left( {0, {R_n} } \right).\] This in conjunction with (\ref{bn}) implies that
\begin{equation}\label{sx2}
\hat \omega _D \left( {R_n } \right) \ge \frac{2}{\pi }e^{ - d_D \left( {0,R_n } \right)}. 
\end{equation}

\noindent
Since ${\Gamma _{{R_n}}}$ denotes the hyperbolic geodesic joining the endpoints of $F_{{R _n}}$ in $\mathbb{D}$, by Theorem \ref{geod} and \cite[p.\ 370]{Beu},
\begin{equation}\nonumber
{\omega _\mathbb{D}}\left( {0,{E_{{R _n}}}} \right) = \frac{1}{2}{\omega _\mathbb{D}}\left( {0,{\Gamma _{{R _n}}}} \right) \le \frac{2}{\pi }{e^{ - {d_\mathbb{D}}\left( {0,{\Gamma _{{R _n}}}} \right)}} 
\end{equation}
and thus
\begin{equation}\label{sx3}
\omega _D \left( {R_n } \right) = \omega _D \left( {0,\psi \left( {E_{R_n } } \right)} \right) \le \frac{2}{\pi }e^{ - d_D \left( {0,\psi \left( {\Gamma _{R_n } } \right)} \right)}.  
\end{equation}
Since  $D$ is symmetric with respect to the real axis, we deduce that 
\[{d_{D}}\left( {{0},\psi \left( {{\Gamma _{{R _n}}}} \right)} \right) = {d_{D}}\left( {0,{r_n}} \right),\]
where ${r_n} = \psi \left( {{\Gamma _{{R_n}}}} \right) \cap \mathbb{R} \in \left( {e^n,e^{n + 1}} \right)$ (see Fig.\ \ref{eik}) and hence by (\ref{sx3}) we conclude that
\begin{equation}\label{sx4}
\omega _D \left( {R_n } \right) \le \frac{2}{\pi }e^{ - d_D \left( 0, r_n \right)}.
\end{equation}

Since $0,\,R_n$ and $r_n$ lie, in this order, along a hyperbolic geodesic (for more details see \cite{Ka}), we have that
\[{d_{D}}\left( {0,{r_n}} \right) = {d_{D}}\left( {0,{R_n}} \right) + {d_{D}}\left( {{R_n},{r_n}} \right)\]
(see \cite[p.\ 14]{Bea}).\ Combining this with (\ref{sx2}) and (\ref{sx4}), we deduce that
\begin{equation}\label{sx5}
\frac{{\hat \omega _D \left( {R_n } \right)}}{{\omega _D \left( {R_n } \right)}} \ge e^{d_D \left( {0,r_n } \right) - d_D \left( {0,R_n } \right)}  = e^{d_D \left( {R_n ,r_n } \right)}.
\end{equation}

\begin{figure}[H] 
	\begin{center}
		\includegraphics[scale=0.6]{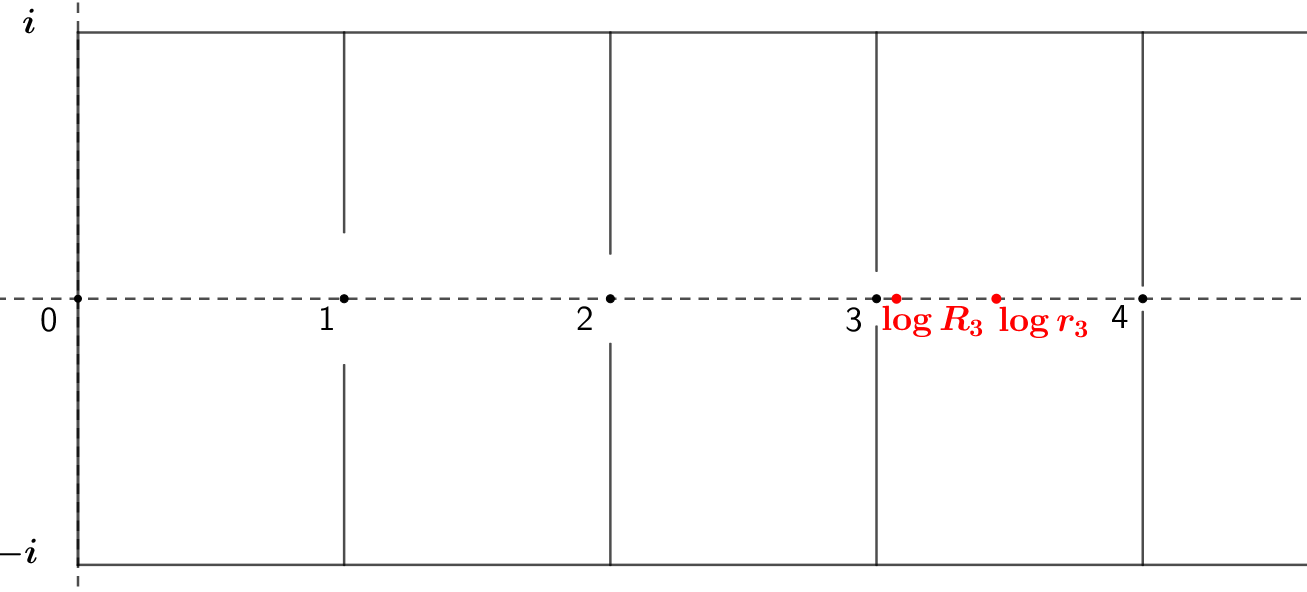}
		\caption{The domain $D'$ and the points $\log R_n$, $\log r_n$ in case $n=3$.}
		\label{re}
	\end{center}
\end{figure}

Now notice that the quasi-hyperbolic distance (see Section \ref{section2}) $\delta _D \left( {R_n ,r_n } \right)$ is equal to $\delta _{D\backslash{\overline{\mathbb{D}} }} \left( {R_n ,r_n } \right)$ because the quasi-hyperbolic geodesic joining $R_n$ to $r_n$ in $D$ and the quasi-hyperbolic geodesic joining $R_n$ to $r_n$ in $D\backslash{\overline{\mathbb{D}} }$ is the segment $\left[ {R_n ,r_n} \right]$ in both cases.\ So, we deduce that
\begin{equation}\label{neo}
d_D \left( {R_n ,r_n } \right) \ge \frac{1}{2}\delta _D \left( {R_n ,r_n } \right) = \frac{1}{2}\delta _{D\backslash{\overline{\mathbb{D}} }} \left( {R_n ,r_n } \right) \ge \frac{1}{4}d_{D\backslash{\overline{\mathbb{D}} }} \left( {R_n ,r_n } \right).
\end{equation}
In order to simplify our computations we use the conformal mapping  $g\left( z \right) = \LOG z$ that maps $D\backslash \overline{\mathbb{D}} $ onto $g\left( {D\backslash \overline{\mathbb{D}} } \right) := D'$ (see Fig.\ \ref{re}).

\noindent
Thus, we get
\begin{eqnarray}\label{anisot}
d_{D\backslash{\overline{\mathbb{D}} }} \left( {R_n ,r_n } \right)&=&{d_{D'}}\left( {\log R_n ,\log r_n } \right) \ge \frac{1}{2}{\delta _{D'}}\left( \log R_n ,\log r_n \right)    \nonumber \\
&=& \frac{1}{2}\int_{{{\mathop{\log R_n}} }}^{\log r_n} {\frac{{dx}}{{d\left( {x,\partial D'} \right)}}}     \ge \frac{1}{2}\int_{\log R_n}^{\log r_n} {\frac{{dx}}{{\sqrt {{{\left( {\frac{1}{{{40^{n}}}}} \right)}^2} + {{\left( {x-n} \right)}^2}} }}} \nonumber \\
&=&\frac{1}{2}\arcsinh \left( {{40^{n}}\left( {\log r_n-n} \right)} \right)-\frac{1}{2}\arcsinh \left( {1} \right) \nonumber \\
&\ge& \frac{1}{2}\arcsinh \left( {{40^{n}} k } \right)-\frac{1}{2}\arcsinh \left( {1} \right), 
\end{eqnarray}
where $k>0$ is a constant independent of $n$ (see \cite{Ka}).\ Now, taking limits in (\ref{anisot}) as $n \to  + \infty $, we obtain
\begin{equation}\nonumber
\mathop {\lim }\limits_{n \to  + \infty } d_{D\backslash{\overline{\mathbb{D}} }} \left( {R_n ,r_n } \right) =  + \infty.
\end{equation}
Thus, by (\ref{neo}) we conclude that 
\begin{equation}\label{sx55}
\mathop {\lim }\limits_{n \to  + \infty } d_{D} \left( {R_n ,r_n } \right) =  + \infty,
\end{equation}
which proves property (ii).\ Finally, by (\ref{sx5}) and (\ref{sx55}), we infer that
\[\mathop {\lim }\limits_{n \to  + \infty } \frac{{{{\hat \omega }_D}\left( {{R_n}} \right)}}{{{\omega _D}\left( {{R_n}} \right)}} =  + \infty \]
and hence property (iii) holds.\ So, there does not exist a positive constant $C$ such that for every $R>0$,
\[{\omega _D}\left( {R} \right) \ge C{{\hat \omega }_D}\left( {R} \right).\]
\qed

\section{Second counter-example}\label{section4}

\begin{figure}[H] 
	\begin{center}
		\includegraphics[scale=0.55]{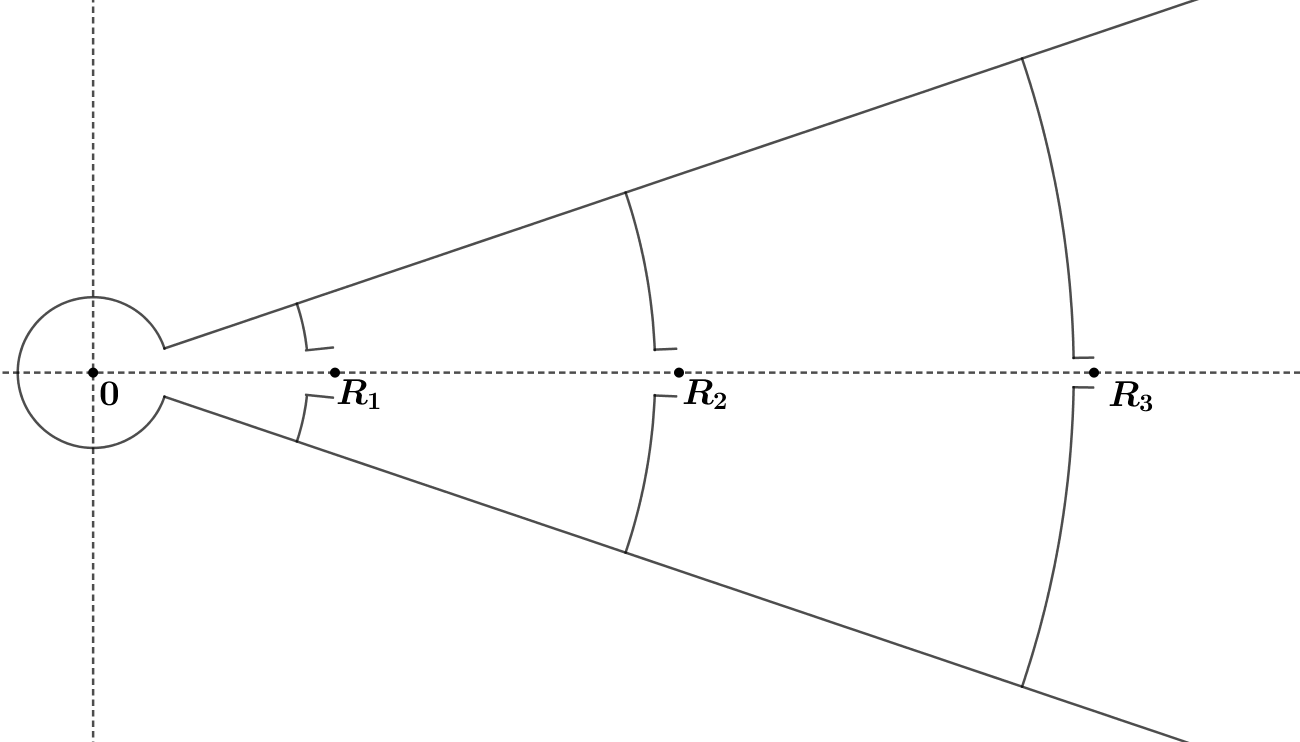}
		\caption{The simply connected domain $D$.}
		\label{dom2}
	\end{center}
\end{figure}

Let $D$ be the simply connected domain of Fig.\ \ref{dom2}, namely,
\[D=\mathbb{D}\cup \left( \left\{ {z \in \mathbb{C}:\left| {\Arg{z}} \right| < 1} \right\}\backslash {D_0}\right),\]
where
\begin{align*}
D_0&= \bigcup_{n=1}^{+\infty} \bigg( \left\{z\in \partial D(0,e^n): \frac{1}{40^n}\leq |\Arg z|\leq 1 \right\}  \\
&\quad\quad\quad\qquad \cup \left\{re^{i\theta}: e^n\leq r\leq e^{n+1/40^n},\, |\theta|=\frac{1}{40^n} \right\}
\bigg).
\end{align*}
We consider the sequence ${\left\{ {{R_n}} \right\}_{n \in \mathbb{N}}}$ with ${R _n} = {e^{n+1/40^n}}$ for every $n \in \mathbb{N}$. 

\begin{theorem}\label{ant2} With the notation above, the simply connected domain $D$ has the following properties:
	\begin{enumerate}[\rm(i)]
		\item $D$ does not satisfy Condition {\rm (2)}.
		\item $D$ satisfies Condition {\rm (1)}.	
		\item $\mathop {\lim }\limits_{n \to  + \infty } \frac{{{{\hat \omega }_D}\left( {{R_n}} \right)}}{{{\omega _D}\left( {{R_n}} \right)}} = + \infty.$	
	\end{enumerate}	
\end{theorem}

\proof  Property (i) is immediate by the construction of $D$.\ So, we prove properties (ii) and (iii).\ First we introduce some notation.\ For $n \in \mathbb{N}$, let ${F_{R_n}}$ be the component of $D \cap \left\{ {\left| z \right| = R_n } \right\}$ that intersects the real axis and $\Gamma_{R_n}$ be the hyperbolic geodesic joining the endpoints of $F _{R _n}$ in $D$.\ Also, we set ${E_{R_n} } = \partial D \cap \left\{ {\left| z \right| \ge R_n } \right\}$ for every $n \in \mathbb{N}$.

\begin{figure}[H] 
	\begin{center}
		\includegraphics[scale=0.65]{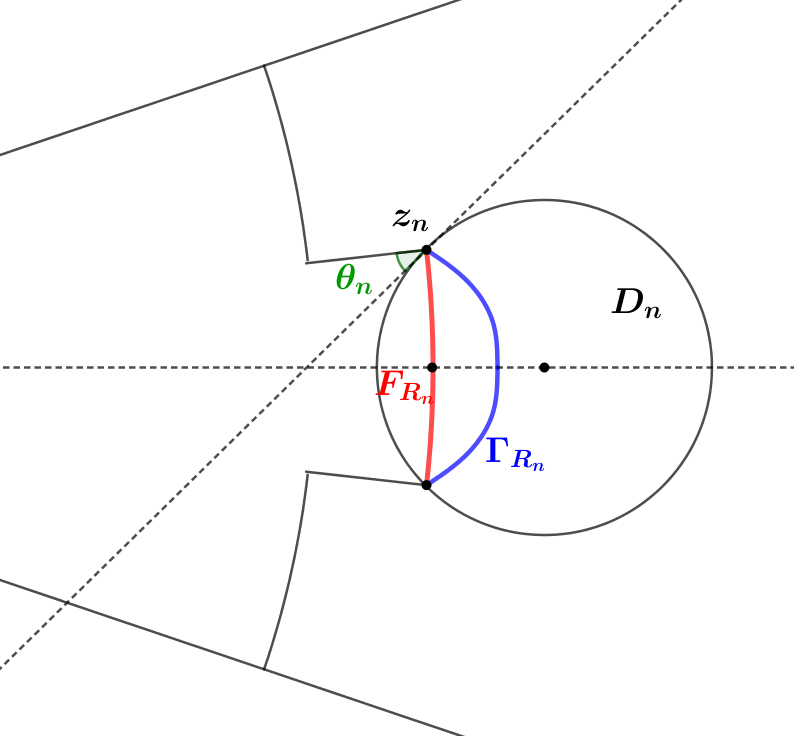}
		\caption{}
		\label{jj}
	\end{center}
\end{figure}

Now we apply J\o rgensen's theorem \cite[p.\ 116]{Jo} that a Euclidean disk inside a simply connected domain is hyperbolically convex.\ Combining this with the construction of $D$, we deduce that, for every $n \in \mathbb{N}$, we can find a disk $D_n \subset D$ centered at a point of $\mathbb{R}$ (see Fig.\ \ref{jj}) that satisfies the following properties:
\begin{enumerate}
	\item $D_n$ contains the arc ${F_{R_n}}$ and the geodesic $\Gamma_{R_n}$.
	\item The endpoints of ${F_{R_n}}$ lie on $\partial D_n$.
	\item  The Euclidean distance of each point of $\partial D_n$ from $\partial D$ is attained on the set $\left\{re^{i\theta}: e^n\leq r\leq e^{n+1/40^n},\, |\theta|=\frac{1}{40^n} \right\}$.
	\item If $\theta_n$ is the acute angle between  $\left\{re^{i\theta}: e^n\leq r\leq e^{n+1/40^n},\, \theta=\frac{1}{40^n} \right\}$ and the tangent of $\partial D_n$ at the point $z_n=e^{n+1/40^n}e^{1/40^ni}$, then $\theta _n  \ge k$ for some constant $k>0$ independent of $n$ (see Fig.\ \ref{jj}).
\end{enumerate}
So, if $s\in D_n \cap \left\{ {z:{\mathop{\rm Im}\nolimits} z \ge 0} \right\}$ then we can easily infer that
\begin{equation}\label{2sx}
\dist\left( {s,\partial D} \right) = \left| {s - z_n } \right|\,\,{\rm or}\,\,\dist\left( {s,\partial D} \right)\ge \sin k\left| {s - z_n } \right|.
\end{equation}
Since $D$ and $D_n$ are symmetric with respect to $\mathbb{R}$, (\ref{2sx}) also holds for every $s\in D_n\cap \left\{ {z:{\mathop{\rm Im}\nolimits} z < 0} \right\}$ by replacing $z_n$ with ${\bar z_n }$.\
So, (\ref{2sx}) in combination with the fact that $F_{R_n},\Gamma_{R_n}$ lie in $D_n$ and join $z_n$ to ${\bar z_n }$ implies that, for every $n \in \mathbb{N}$, the quasi-hyperbolic distance between any point of $F_{R_n}$ and $\Gamma_{R_n}$ is bounded from above by an absolute positive constant.\ Thus, the hyperbolic distance between any point of $F_{R_n}$ and $\Gamma_{R_n}$ is bounded from above by an absolute positive constant.\ This proves property (ii). 

\begin{figure}[H] 
	\begin{center}
		\includegraphics[scale=0.5]{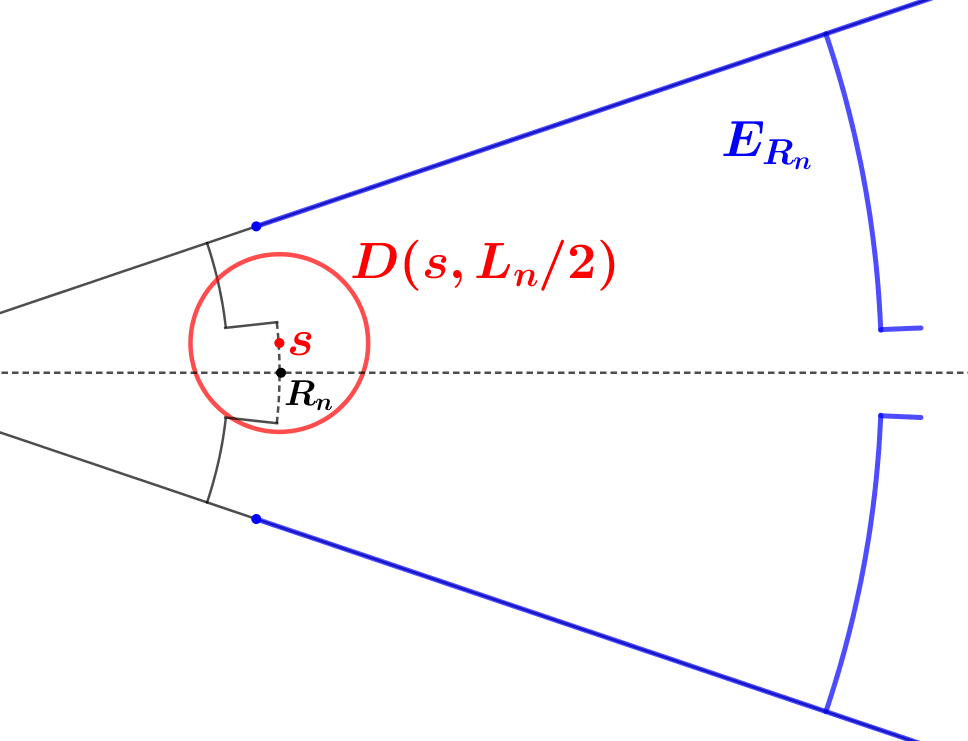}
		\caption{}
		\label{antip22}
	\end{center}
\end{figure}

Now set $L_n  = \dist \left( {R_n ,\left\{ {z \in \mathbb{C}:\left| {\Arg z} \right| = 1} \right\}} \right)$ and $d_n  = \dist\left( {R_n ,\partial D} \right)$.\ By the construction of $D$, there exists a number $n_0 \in \mathbb{N}$ such that for every $n>n_0$ and every $s \in F_{R_n}$ (see Fig.\ \ref{antip22}), 
\[D\left( {s,\frac{{L_n }}{2}} \right) \subset \left\{ {z \in \mathbb{C}:\left| {\Arg z} \right| < 1} \right\}\,\,\,{\rm and}\,\,\,D\left( {s,\frac{{L_n }}{2}} \right) \cap E_{R_n }  = \emptyset.\]
Fix a number $n>n_0$ and a point $s \in F_{R_n}$.\ The Riemann mapping theorem implies that there exists a conformal mapping $f$ from $\mathbb{D}$ onto $D$ such that $f\left( 0 \right) = s$.\ Therefore, applying Theorem \ref{bonk} with its notation, we have that
\begin{eqnarray}\label{e2}
\omega _D \left( {s,E_{R_n } } \right) &=& \omega _\mathbb{D} \left( {0,f^{ - 1} \left( {E_{R_n } } \right)} \right) \le {\rm cap}f^{ - 1} \left( {E_{R_n } } \right)\le {\rm cap}E_f \left( {\frac{{L_n }}{2}} \right)  \nonumber \\
&\le& K\sqrt {\frac{{2d_s }}{{L_n }}}  \le K\sqrt {\frac{{2d_n }}{{L_n }}}, 
\end{eqnarray}
where $d_s  = \dist\left( {s,\partial D} \right)$.\ So, by Theorem \ref{markov} and relation (\ref{e2}), we infer that for every $n>n_0$,
\[\frac{{\hat \omega _D \left( {R_n } \right)}}{{\omega _D \left( {R_n } \right)}} = \frac{{\omega _D \left( {0,F_{R_n } } \right)}}{{\omega _D \left( {0,E_{R_n } } \right)}} = \frac{{\omega _D \left( {0,F_{R_n } } \right)}}{{\int_{F_{R_n } } {\omega _D \left( {0,ds} \right)\omega _D \left( {s,E_{R_n } } \right)} }} \ge \frac{1}{K}\sqrt {\frac{{L_n }}{{2d_n }}}. \]
\noindent
Taking limits as $n\to +\infty$, we deduce that
\[\mathop {\lim }\limits_{n \to  + \infty } \frac{1}{K}\sqrt {\frac{{L_n }}{{2d_n }}} =  + \infty\]
and hence
\[\mathop {\lim }\limits_{n \to  + \infty } \frac{{\hat \omega _D \left( {R_n } \right)}}{{\omega _D \left( {R_n } \right)}} =  + \infty.\]
This proves property (iii).

\qed

\section{Proof of theorem \ref{kyrio}}\label{section5}

\proof[Proof of Theorem \ref{kyrio}] 

\begin{figure}[H]
	\begin{minipage}{0.5\textwidth}
		\begin{center}
			\includegraphics[scale=0.55]{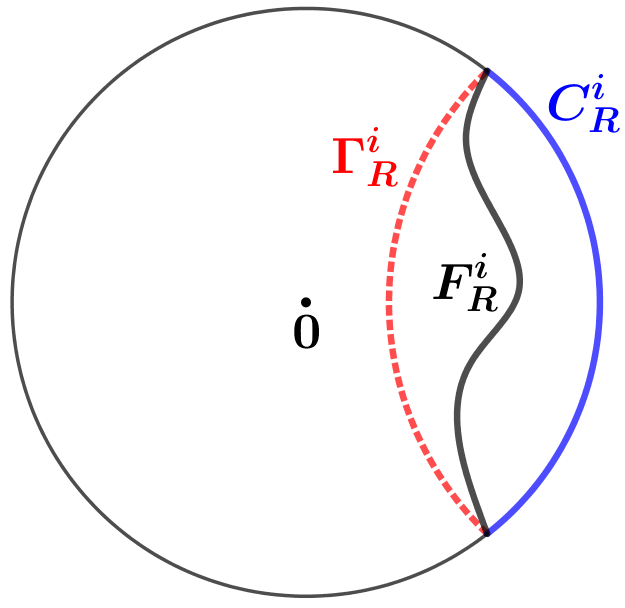}
			\vspace*{0.45cm}
			\caption{}
			\label{int1}
		\end{center}
	\end{minipage}\hfill
	\begin{minipage}{0.5\textwidth}
		\begin{center}
			\includegraphics[scale=0.55]{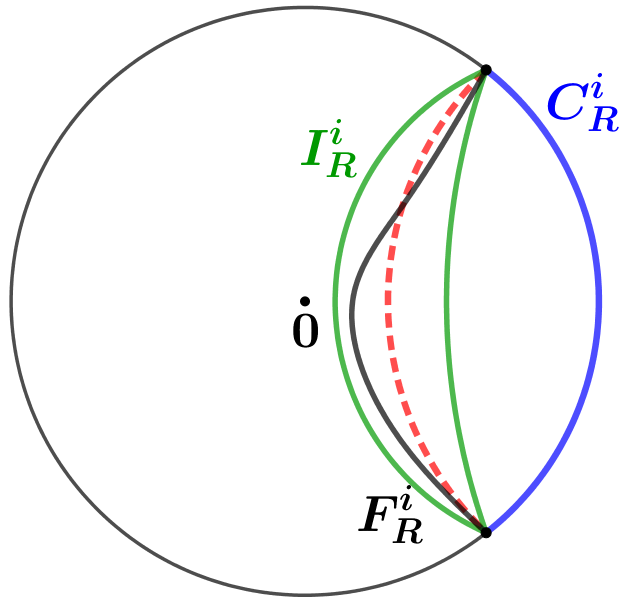}
			\vspace*{0.3cm}
			\caption{}
			\label{int2}
		\end{center}
	\end{minipage}
\end{figure}

Since $D$ is a simply connected domain, the Riemann mapping theorem implies that there exists a conformal mapping $\psi$ from $\mathbb{D}$ onto $D$ with $\psi \left( 0 \right) = 0$.\ Now we introduce some notation.\ For $R >0$, we set ${F_R } = \left\{ {z \in \mathbb{D}:\left| {\psi \left( z \right)} \right| = R } \right\}$, that is, $\psi \left( {{F_R}} \right) = D \cap \left\{ {z:\left| z \right| = R} \right\}$.\ Note that since $\psi \left( {{F_R}} \right)$ is a countable union of open arcs in $D$ that are the intersection of $D$ with the circle $\left\{ {z:\left| z \right| = R} \right\}$, the preimage of every such arc is also an arc in $\mathbb{D}$ with two distinct endpoints on $\partial \mathbb{D}$ (see Proposition 2.14 \cite[p.\ 29]{Pom}).\ Also, let $N\left( R \right) \in \mathbb{N} \cup \left\{ { + \infty } \right\}$ denote the number of components of $F_R$ and

\[I_R = \left\{ \begin{array}{l}
\left\{ {1,2, \ldots ,N\left( R \right)} \right\},\,{\rm{if}}\,N\left( R \right) <  + \infty  \\ 
\mathbb{N},\,\,\,\,\,\,\,\,\,\,\,\,\,\,\,\,\,\,\,\,\,\,\,\,\,\,\,\,\,\,\,\,\,\,\,\,\,\,\,\,\,{\rm{if}}\,N\left( R \right) =  + \infty  \\ 
\end{array} \right..\]
If ${\left\{ {F_R^i } \right\}_{i \in I_R } }$ are the components of $F_R$, then, for every $i \in I_R$, we set ${\Gamma}_R^i$ be the hyperbolic geodesic joining the endpoints of $F_R^i $ in $\mathbb{D}$ and $C_R^i$ be the arc of $\partial \mathbb{D}$ joining the endpoints of ${\Gamma}_R^i$ and lying on the boundary of the component of $\mathbb{D}\backslash {{\Gamma}_R^i}$ which does not contain the origin (see Fig.\ \ref{int1}).

Suppose that Conditions (2) and (3) are satisfied.\ For every $R>0$ and $i \in I_R$ and for each $z\in {{\Gamma}_R^i}$, we have that 
\begin{equation}\nonumber
{\omega _\mathbb{D}}\left( {z,C_R^i} \right)=\frac{1}{2},
\end{equation} 
(see \cite[p.\ 370]{Beu}).\ Condition (3) implies that each crosscut $F_R^i $ is contained in the component of $\mathbb{D}\backslash {{\Gamma}_R^i}$ bounded by ${{\Gamma}_R^i}$ and $C_R^i$ (see Fig.\ \ref{int1}).\ Thus, by the maximum principle, we deduce that for every $z\in {F_R^i}$,
\begin{equation}\label{a1}
{\omega _\mathbb{D}}\left( {z,C_R^i} \right)\ge\frac{1}{2}.
\end{equation}
\noindent
Applying Theorem \ref{markov} and relation (\ref{a1}), we infer that
\begin{eqnarray}\label{a22}
\omega _\mathbb{D} \left( {0,C_R^i } \right) &=& \int\limits_{F_R^i } {\omega _\mathbb{D} \left( {0,dz} \right)\omega _\mathbb{D} \left( {z,C_R^i } \right)}  \ge \frac{1}{2}\int\limits_{F_R^i } {\omega _\mathbb{D} \left( {0,dz} \right)} \nonumber \\ 
&=& \frac{1}{2}\omega _\mathbb{D} \left( {0,F_R^i } \right)
\end{eqnarray}
for every $R>0$ and every $i \in I_R$.\ Condition (2) and the conformal invariance of harmonic measure imply that 
\[
\omega _D \left( R \right) \ge  \omega _\mathbb{D} \left( {0,\bigcup\limits_{i \in I_R } {C_R^i } } \right)=
\sum\limits_{i \in I_R } {\omega _\mathbb{D} \left( {0,C_R^i } \right)}.\]
Combining this with (\ref{a22}) we get
\begin{eqnarray}
\omega _D \left( R \right) &\ge& \sum\limits_{i \in I_R } {\omega _\mathbb{D} \left( {0,C_R^i } \right)}  \ge \frac{1}{2}\sum\limits_{i \in I_R } {\omega _\mathbb{D} \left( {0,F_R^i } \right)}= \frac{1}{2}\omega _\mathbb{D} \left( {0, \bigcup\limits_{i \in I_R } {F_R^i }} \right) \nonumber \\ 
&=& \frac{1}{2}\omega _\mathbb{D} \left( {0, {F_R }} \right)  = \frac{1}{2}\omega _D \left( {0,\psi \left( {F_R } \right)} \right) = \frac{1}{2}\hat \omega _D \left( R \right) \nonumber
\end{eqnarray}
and thus we have the desired result
\[{{\hat \omega }_D}\left( {R} \right) \le 2{\omega _D}\left( {R} \right)\]
for every $R>0$. 

Now suppose that Conditions (1) and (2) are satisfied.\ By Condition (1) we infer that, for every $R>0$ and every $i \in I_R$, there exists a hyperbolic $k_c$-neighborhood, $U_R^i$, of ${{\Gamma}_R^i}$ such that $\partial{U_R^i}$ consists of two circular arcs in $\mathbb{D}$ and ${F_R^i }$ is contained in $U_R^i$ (see Fig.\ \ref{int2}).\ Note that $k_c$ is a positive constant that depends only on $c$.\ Let ${I_R^i }$ denote the circular arc of $\partial{U_R^i}$ such that ${F_R^i }$ is contained in the component of $\mathbb{D}\backslash {I_R^i}$ bounded by ${I_R^i}$ and $C_R^i$ (see Fig.\ \ref{int2}).\ For every $R>0$ and $i \in I_R$ and for each $z\in {I_R^i}$, we have that 
\begin{equation}\nonumber
{\omega _\mathbb{D}}\left( {z,C_R^i} \right)=k',
\end{equation}
where $k'$ lies in the open interval $\left(0,1 \right)$ and depends only on $k_c$ and hence only on $c$.\ Now we repeat the argument above letting ${I_R^i}$ play the role of ${{\Gamma}_R^i}$.\ Therefore, for every $z\in {F_R^i}$,
\begin{equation}\label{a11}
{\omega _\mathbb{D}}\left( {z,C_R^i} \right)\ge k'.
\end{equation}
By Theorem \ref{markov} and relation (\ref{a11}), we infer that
\begin{equation}\nonumber
\omega _\mathbb{D} \left( {0,C_R^i } \right) = \int\limits_{F_R^i } {\omega _\mathbb{D} \left( {0,dz} \right)\omega _\mathbb{D} \left( {z,C_R^i } \right)}  \ge k'\omega _\mathbb{D} \left( {0,F_R^i } \right)
\end{equation}
for every $R>0$ and every $i \in I_R$.\ This in conjunction with Condition (2) implies that 
\begin{eqnarray}
\omega _D \left( R \right) &\ge& \omega _\mathbb{D} \left( {0,\bigcup\limits_{i \in I_R } {C_R^i } } \right)=
\sum\limits_{i \in I_R } {\omega _\mathbb{D} \left( {0,C_R^i } \right)}\ge k' \sum\limits_{i \in I_R } {\omega _\mathbb{D} \left( {0,F_R^i } \right)} \nonumber \\ 
&=& k'\omega _\mathbb{D} \left( {0, {F_R }} \right)= k'\hat \omega _D \left( R \right).\nonumber
\end{eqnarray}
So, we conclude that for every $R>0$,
\[{{\hat \omega }_D}\left( {R} \right) \le K{\omega _D}\left( {R} \right),\]
where $K = \frac{1}{{k'}}$ is a positive constant that depends only on $c$. 
      
\qed

\section{Proof of theorem \ref{star}}\label{section6}

In the proof of Theorem \ref{star} we will use the following result which is an easy computation coming from the conformal invariance of harmonic measure.

\begin{lemma}\label{le2} Let $a\in \left( {0,1} \right)$ and $b\in \left[ {0,1} \right)$.\ Then
	\[{\omega _{\mathbb{D} \backslash \left[ {a,1} \right)}}\left( { - b,\partial {\mathbb{D}}} \right) = 1 - \frac{2}{\pi }\arctan \frac{1}{{\sqrt {{{\left( {\frac{{\left( {1 + a} \right)\left( {1 + b} \right)}}{{\left( {1 - a} \right)\left( {1 - b} \right)}}} \right)}^2} - 1} }}.\]
\end{lemma}

\begin{figure}[H] 
	\begin{center}
		\includegraphics[scale=0.65]{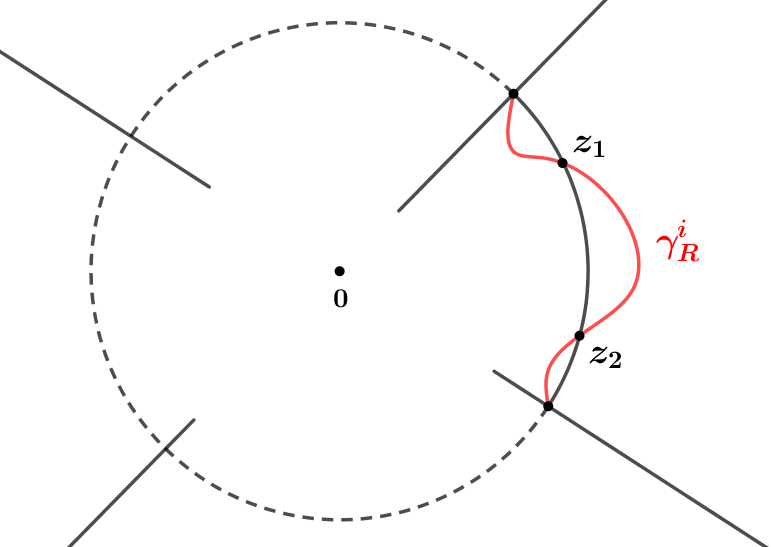}
		\caption{}
		\label{eiko}
	\end{center}
\end{figure}

\proof[Proof of Theorem \ref{star}] Let $D$ be a starlike domain in $\mathbb{C}$.\ Using the notation of the proof of Theorem \ref{kyrio}, we will prove that Conditions (2) and (3) are satisfied.\ Since $D$ is starlike, Condition (2) is obviously satisfied and thus we prove  Condition (3).\ Let ${F_R^i}$ be a component of $F_R$ for some $i \in I_R$.\ Suppose that $\psi \left( {\Gamma _R^i} \right) \not\subset \overline D  \cap \left\{ {z:\left| z \right| \le R} \right\}$, then $\psi \left( {\Gamma _R^i} \right)$ contains a curve ${\gamma _R^i}$ lying in $D\backslash D\left( {0,R} \right)$ with endpoints $z_1,z_2 \in \partial {D \left( {0,R } \right)}$ (see Fig.\ \ref{eiko}).\ Since $\psi \left( {\Gamma _R^i} \right)$ is the hyperbolic geodesic joining the endpoints of $\psi \left( {F_R^i} \right)$ in $D$, ${\gamma _R^i}$ is the hyperbolic geodesic joining $z_1$ to ${z_2}$ in $D$.\ Notice that $D$ is a hyperbolic region in $\mathbb{C}$ such that $D \cap \partial D\left( {0,R} \right) \ne \emptyset $.\ Since $D$ is starlike, we have that $D \backslash D\left( {0,R } \right) \subset D ^ *$, where $D^ *$ is the reflection of $D$ in the circle $\partial D\left( {0,R} \right)$.\ So, applying Theorem \ref{rp}, we get
\[{\lambda _{{D ^ * }}}\left( z \right) < {\lambda _D }\left( z \right),\;z \in {\gamma _R^i}\]
and thus
\[\int_{{{\gamma _R^i}}^ * } {{\lambda _D }\left( {{z^ * }} \right)\left| {d{z^ * }} \right|}  < \int_{{\gamma _R^i}} {{\lambda _D }\left( z \right)\left| {dz} \right|}, \]
where ${{\gamma _R^i}}^ * $ is the reflection of ${\gamma _R^i}$ in $\partial D\left( {0,R} \right)$.\ But this leads to contradiction because ${\gamma _R^i}$ is the hyperbolic geodesic joining $z_1$ to ${z_2}$ in $D$.\ So, $\psi \left( {\Gamma _R^i} \right) \subset \overline D  \cap \left\{ {z:\left| z \right| \le R} \right\}$ and thus Condition (3) is satisfied.\ Theorem \ref{kyrio} implies that for every $R>0$,
\[{{\hat \omega }_D}\left( {R} \right) \le 2{\omega _D}\left( {R} \right).\]
Now we prove that the constant $2$ is best possible.\ Consider the Koebe function $K \left( z \right) = \frac{z}{{{{\left( {1 - z} \right)}^2}}}$ which maps $\mathbb{D}$  conformally onto $D_0:=\mathbb{C}\backslash \left( { - \infty , - \frac{1}{4}} \right]$.\ For $R>{\frac{1}{4}}$, by the conformal invariance of harmonic measure and Lemma \ref{le2}, we have

\begin{eqnarray}\label{ko1}
{{\hat \omega }_{D_0}}\left( {R} \right)&=& {\omega _{D\left( {0,R} \right)\backslash \left( { - R, - \frac{1}{4}} \right]}}\left( {0,\partial D\left( {0,R} \right)} \right)={\omega _{\mathbb{D}\backslash \left( { - 1, - \frac{1}{{4R}}} \right]}}\left( {0,\partial \mathbb{D}} \right) \nonumber \\
&=&{\omega _{\mathbb{D}\backslash \left[ {\frac{1}{{4R}},1} \right)}}\left( {0,\partial \mathbb{D}} \right)= 1 - \frac{2}{\pi }\arctan \frac{1}{{\sqrt {{{\left( {\frac{{4R + 1}}{{4R - 1}}} \right)}^2} - 1} }} \nonumber \\
&=& 1 - \frac{2}{\pi }\arctan \frac{{4R - 1}}{{4\sqrt R }}. 
\end{eqnarray}
Using the fact that
\[{K^{ - 1}}\left( { - R} \right) = \left( {\frac{{2R - 1}}{{2R}}} \right) \pm i\frac{{\sqrt {4R - 1} }}{{2R}}\]
and the conformal invariance of harmonic measure, we deduce that
\begin{eqnarray}\label{ko2}
{{\omega }_{D_0}}\left( {R} \right)&=& {\omega _{D_0}}\left( {0,\left( { - \infty , - R} \right]} \right) \nonumber \\
&=&{\omega _\mathbb{D}}\left( {0,\arc\left( {\left( {\frac{{2R - 1}}{{2R}}} \right) - i\frac{{\sqrt {4R - 1} }}{{2R}},\left( {\frac{{2R - 1}}{{2R}}} \right) + i\frac{{\sqrt {4R - 1} }}{{2R}}} \right)} \right) \nonumber \\
&=&2{\omega _\mathbb{D}}\left( {0,\arc\left( {1,\left( {\frac{{2R - 1}}{{2R}}} \right) + i\frac{{\sqrt {4R - 1} }}{{2R}}} \right)} \right) \nonumber \\
&=&\frac{1}{\pi }\arctan \frac{{\sqrt {4R - 1} }}{{2R - 1}},
\end{eqnarray}
where $\arc\left( {\left( {\frac{{2R - 1}}{{2R}}} \right) - i\frac{{\sqrt {4R - 1} }}{{2R}},\left( {\frac{{2R - 1}}{{2R}}} \right) + i\frac{{\sqrt {4R - 1} }}{{2R}}} \right)$ denotes the arc of $\partial \mathbb{D}$ joining $\left( {\frac{{2R - 1}}{{2R}}} \right) - i\frac{{\sqrt {4R - 1} }}{{2R}}$ to $\left( {\frac{{2R - 1}}{{2R}}} \right) + i\frac{{\sqrt {4R - 1} }}{{2R}}$ counterclockwise.\ Applying (\ref{ko1}) and (\ref{ko2}), we infer that
\[\mathop {\lim }\limits_{R \to  + \infty } \frac{{{{\hat \omega }_{D_0}}\left( {R} \right)}}{{{\omega _{D_0}}\left( {R} \right)}} = \mathop {\lim }\limits_{R \to  + \infty } \frac{{\pi  - 2\arctan \frac{{4R - 1}}{{4\sqrt R }}}}{{\arctan \frac{{\sqrt {4R - 1} }}{{2R - 1}}}} = 2.\]
Suppose that there exists a positive constant $C<2$ such that for every starlike domain $D$ and every $R>0$, ${{\hat \omega }_D}\left( {R} \right) \le C{\omega _D}\left( {R} \right)$.\ This implies for $D_0$ that
\[2 = \mathop {\lim }\limits_{R \to  + \infty } \frac{{{{\hat \omega }_{{D_0}}}\left( {R} \right)}}{{{\omega _{{D_0}}}\left( {R} \right)}} \le C,\]
which leads to contradiction.\ Therefore, the constant $2$ is best possible.

\qed
\\

Note that we could also prove Theorem \ref{star}  by using instead of Minda's reflection principle and Theorem \ref{kyrio}, the strong Markov property for harmonic measure (see Section \ref{section2}).

\begin{figure}[H]
		\begin{center}
			\includegraphics[scale=0.7]{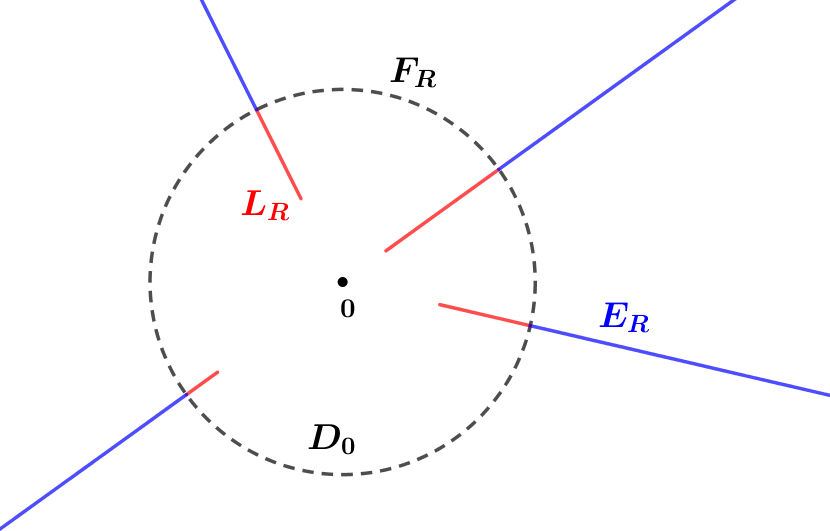}
			\caption{The starlike domain $D$.}
			\label{mark}
		\end{center}
\end{figure}

\proof[Another proof of Theorem \ref{star}] Let $D$ be a starlike domain in $\mathbb{C}$.\ Set ${F_R} = D \cap \partial D\left( {0,R} \right)$, ${E_R} = \partial D\backslash D\left( {0,R} \right)$, ${L_R} = \partial D \cap D\left( {0,R} \right)$ and ${D_0} = D \cap D\left( {0,R} \right)$ as illustrated in Fig.\ \ref{mark}.\ So, we have the relations
\begin{equation}\label{ss1}
{{\hat \omega }_D}\left( {R} \right) = {\omega _{{D_0}}}\left( {0,{F_R}} \right) = 1 - {\omega _{{D_0}}}\left( {0,{L_R}} \right)
\end{equation}
and
\begin{equation}\label{ss2}
{\omega _D}\left( {R} \right) = {\omega _D}\left( {0,{E_R}} \right) = 1 - {\omega _D}\left( {0,{L_R}} \right).
\end{equation}
\noindent
By Theorem \ref{markov}, 
\[{\omega _D}\left( {0,{L_R}} \right) = {\omega _{{D_0}}}\left( {0,{L_R}} \right) + \int_{{F_R}}{{\omega _{{D_0}}}\left( {0,ds} \right){\omega _D}\left( {s,{L_R}} \right)}\]
which in conjunction with (\ref{ss1}) and (\ref{ss2}) implies that
\begin{equation}\label{ss3}
{{\hat \omega }_D}\left( {R} \right) = {\omega _D}\left( {R} \right) + \int_{{F_R}}{{\omega _{{D_0}}}\left( {0,ds} \right){\omega _D}\left( {s,{L_R}} \right)}.
\end{equation}
Let $N\left( R \right) \in \mathbb{N} \cup \left\{ { + \infty } \right\}$ denote the number of components of ${F_R}$ and 

\[I_R = \left\{ \begin{array}{l}
\left\{ {1,2, \ldots ,N\left( R \right)} \right\},\,{\rm{if}}\,N\left( R \right) <  + \infty  \\ 
\mathbb{N},\,\,\,\,\,\,\,\,\,\,\,\,\,\,\,\,\,\,\,\,\,\,\,\,\,\,\,\,\,\,\,\,\,\,\,\,\,\,\,\,\,{\rm{if}}\,N\left( R \right) =  + \infty  \\ 
\end{array} \right..\]
If ${\left\{ {F_R^i } \right\}_{i \in I_R } }$ are the components of $F_R$, then
\begin{equation}\label{ss4}
\int_{{F_R}} {{\omega _{{D_0}}}\left( {0,ds} \right){\omega _D}\left( {s,{L_R}} \right)}  = \sum\limits_{i \in I_R} {\int_{F_R^i} {{\omega _{{D_0}}}\left( {0,ds} \right){\omega _D}\left( {s,{L_R}} \right)} },
\end{equation}
since ${F_R^i}$ are mutually disjoint sets.

\begin{figure}[H]
	\begin{center}
		\includegraphics[scale=0.65]{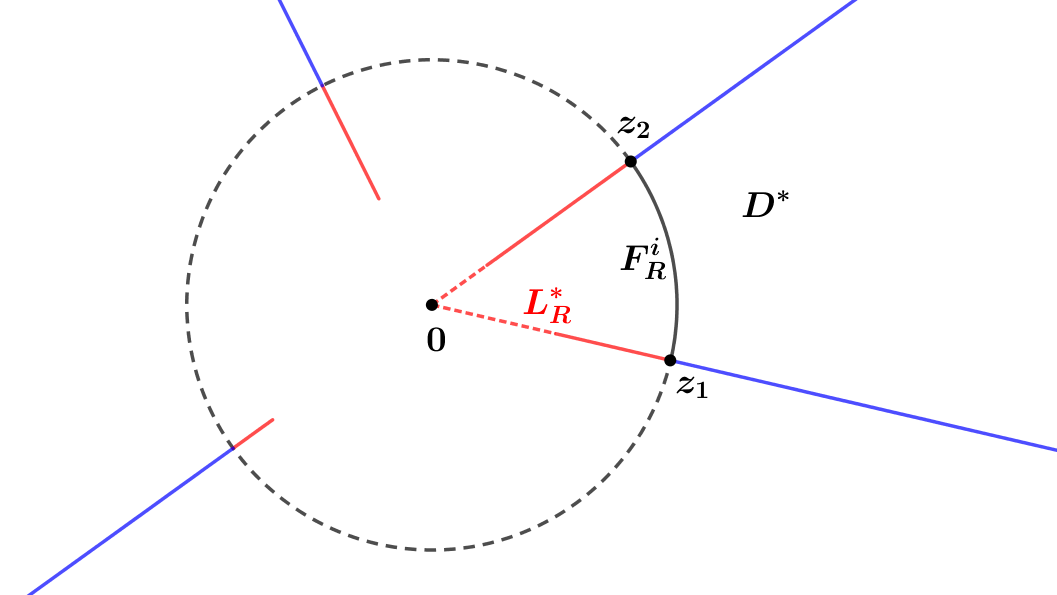}
		\caption{}
		\label{mark2}
	\end{center}
\end{figure}

Now let ${F_R^i}$ be a component of $F_R$ for some $i \in I_R$.\ If $z_1,z_2$ denote the endpoints of ${F_R^i}$ such that $\Arg{z_1} < \Arg{z_2}$, then we set 
\[L_R^ *  = \left\{ {r{e^{i\Arg{z_1}}}:0 \le r \le R} \right\} \cup \left\{ {r{e^{i\Arg{z_2}}}:0 \le r \le R} \right\}\]
and
\[{D^ * } = \left\{ {z \in \mathbb{C}:\Arg{z_1} < \Arg z < \Arg{z_2}} \right\}\]
as illustrated in Fig.\ \ref{mark2}.\ For every $s \in F_R^i$,
\begin{equation}\label{ss5}
{\omega _D}\left( {s,{L_R}} \right) \le {\omega _{{D^ * }}}\left( {s,L_R^ * } \right).
\end{equation}

\noindent
If $\theta  = \Arg{z_2} - \Arg{z_1}$, we consider the conformal mappings 
\[{f_1}\left( z \right) = z{e^{ - i\left( {\Arg{z_2} - {\theta  \mathord{\left/
					{\vphantom {\theta  2}} \right.
					\kern-\nulldelimiterspace} 2}} \right)}},\,{f_2}\left( z \right) = {z^{{\pi  \mathord{\left/
				{\vphantom {\pi  \theta }} \right.
				\kern-\nulldelimiterspace} \theta }}},\, 
{f_3}\left( z \right) = \frac{{z - 1}}{{z + 1}}.\] Then the composition $f = {f_3} \circ {f_2} \circ {f_1}$ maps $D^ * $ conformally onto $\mathbb{D}$.\ Since $f\left( {F_R^i} \right)$ is the hyperbolic geodesic joining $f\left( {{z_1}} \right)$ to $f\left( {{z_2}} \right)$ in $\mathbb{D}$, for every $s \in F_R^i$,
\[{\omega _{{D^ * }}}\left( {s,L_R^ * } \right) = {\omega _\mathbb{D}}\left( {f\left( s \right),f\left( {L_R^ * } \right)} \right) = \frac{1}{2}.\]
This in combination with (\ref{ss5}) implies that for every $s \in F_R^i$,

\[{\omega _D}\left( {s,{L_R}} \right) \le \frac{1}{2}.\]
By this and relations (\ref{ss3}) and (\ref{ss4}) we infer that
\begin{eqnarray}
{{\hat \omega }_D}\left( {R} \right) &=& {\omega _D}\left( {R} \right) + \sum\limits_{i \in I_R} {\int_{F_R^i} {{\omega _{{D_0}}}\left( {0,ds} \right){\omega _D}\left( {s,{L_R}} \right)} } \nonumber \\ 
&\le& {\omega _D}\left( {R} \right) + \frac{1}{2}\sum\limits_{i \in I_R} {\int_{F_R^i} {{\omega _{{D_0}}}\left( {0,ds} \right)} } \nonumber \\
&=& {\omega _D}\left( {R} \right)+\frac{1}{2}\sum\limits_{i\in I_R} {\omega _{{D_0}}}\left( {0,F_R^i} \right)= {\omega _D}\left( {R} \right)+\frac{1}{2}{\omega _{{D_0}}}\left( {0,{F_R}} \right) \nonumber \\
&=&{\omega _D}\left( {R} \right)+\frac{1}{2}{{\hat \omega }_D}\left( {R} \right), \nonumber
\end{eqnarray}
and thus for every $R>0$,
\[{{\hat \omega }_D}\left( R \right) \le 2{\omega _D}\left( R \right).\]
The fact that the constant $2$ is best possible is proved as before.

\qed


\begin{bibdiv}
\begin{biblist}

\bib{Ahl}{book}{
title={Conformal Invariants: Topics in Geometric Function Theory},
author={L.V. Ahlfors},
date={1973},
publisher={McGraw-Hill},
address={New York}
}
\bib{Bae}{article}{
	title={The size of the set where a univalent function is large},
	author={A. Baernstein},
	journal={J. d' Anal. Math.},
	volume={70},
	date={1996},
	pages={157--173}
}
\bib{Bo}{article}{
	title={Lengths of radii under conformal maps of the unit disk},
	author={Z. Balogh and M. Bonk},
	journal={Proc. Amer. Math. Soc.},
	volume={127},
	date={1999},
	pages={801--804}
}
\bib{Bea}{article}{
title={The hyperbolic metric and geometric function theory},
author={A.F. Beardon and D. Minda,},
journal={Quasiconformal mappings and their applications},
date={2007},
pages={9--56}
}
\bib{Be}{article}{
	title={Harmonic measure on simply connected domains of fixed inradius},
	author={D. Betsakos},
	journal={Ark. Mat.},
	volume={36},
	date={1998},
	pages={275--306}
}
\bib{Bet}{article}{
	title={Geometric theorems and problems for harmonic measure},
	author={D. Betsakos},
	journal={Rocky Mountain J. of Math.},
	volume={31},
	date={2001},
	pages={773--795}
}
\bib{Beu}{book}{
	title={The Collected Works of Arne Beurling},
	subtitle={Vol. 1, Complex Analysis},
	author={A. Beurling},
	date={1989},
	publisher={Birkh\"{a}user},
	address={Boston}
}
\bib{Es}{article}{
	title={On analytic functions which are in $H^p$ for some positive $p$},
	author={M. Ess{\' e}n},
	journal={Ark. Mat.},
	volume={19},
	date={1981},
	pages={43--51}
}
\bib{Ess}{article}{
	title={Harmonic majorization and classical analysis},
	author={M. Ess{\' e}n K. Haliste J.L. Lewis and D.F. Shea},
	journal={J. London Math. Soc.},
	volume={32},
	date={1985},
	pages={506--520}
}
\bib{Esse}{article}{
	title={Harmonic majorization and thinness },
	author={M. Ess{\' e}n},
	journal={Proc. of the 14th Winter School on Abstract Analysis},
	volume={},
	date={1987},
	pages={295--304}
}
\bib{Gar}{book}{
title={Harmonic Measure},
author={J.B. Garnett and D.E. Marshall},
date={2005},
publisher={Cambridge University Press},
address={Cambridge}
}
\bib{Hay}{book}{
	title={Subharmonic Functions},
	subtitle={Volume 2},
	author={W.K. Hayman},
	date={1989},
	publisher={Academic Press},
	address={London}
}
\bib{We}{article}{
	title={On the coefficients and means of functions omitting values},
	author={W.K. Hayman and A. Weitsman},
	journal={Math. Proc. Cambridge Philos. Soc.},
	volume={77},
	date={1975},
	pages={119--137}
}
\bib{Jo}{article}{
	title={On an inequality for the hyperbolic measure and its applications in the theory of functions},
	author={V. J\o rgensen},
	journal={Math. Scand.},
	volume={4},
	date={1956},
	pages={113--124}
}
\bib{Ka}{article}{
title={On a relation between harmonic measure and hyperbolic distance on planar domains},
author={C. Karafyllia},
journal={Indiana Univ. Math. J. (to appear)}
}
\bib{Kara}{article}{
	title={On the Hardy number of a domain in terms of harmonic measure and hyperbolic distance},
	author={C. Karafyllia},
	journal={(submitted). ArXiv:1908.11845}
}
\bib{Kim}{article}{
title={Hardy spaces and unbounded quasidisks},
author={Y.C. Kim and T. Sugawa},
journal={Ann. Acad. Sci. Fenn.},
volume={36},
date={2011},
pages={291--300}
}
%\bib{Le}{book}{
%	title={Quasiconformal Mappings in the Plane},
%	author={O. Lehto and K.I. Virtanen},
%	date={1973},
%	publisher={Springer-Verlag},
%	address={New York Heidelberg Berlin}
%}
\bib{Mi}{article}{
	title={Inequalities for the hyperbolic metric and applications to geometric function theory},
	author={D. Minda},
	journal={Lecture Notes in Math.},
	volume={1275},
	date={1987},
	pages={235--252}
}
\bib{Co}{article}{
title={Geometric models, iteration and composition operators},
author={P. Poggi-Corradini},
journal={Ph.D. Thesis, University of Washington},
date={1996}
}
\bib{Co1}{article}{
	title={The Hardy class of geometric models and the essential spectral radius of composition operators},
	author={P. Poggi-Corradini},
	journal={Journal of Functional Analysis},
	volume={143},
	date={1997},
	pages={129--156}
}
\bib{Co2}{article}{
	title={The Hardy class of K{\oe}nigs maps},
	author={P. Poggi-Corradini},
	journal={Michigan Math. J.},
	volume={44},
	date={1997},
	pages={495--507}
}
\bib{Pom}{book}{
	title={Boundary Behaviour of Conformal Maps},
	author={C. Pommerenke},
	date={1992},
	publisher={Springer-Verlag},
	address={Berlin}
}
\bib{Port}{book}{
	title={Brownian Motion and Classical Potential Theory},
	author={S.C. Port and C.J. Stone},
	date={1978},
	publisher={Academic Press},
	address={New York}
}
\bib{Ra}{book}{
	title={Potential Theory in the Complex Plane},
	author={T. Ransford},
	date={1995},
	publisher={Cambridge University Press},
	address={Cambridge}
}
\bib{Sa}{article}{
	title={Isoperimetric inequalities for the least harmonic majorant of $\left| x \right|^p$},
	author={M. Sakai},
	journal={Trans. Amer. Math. Soc.},
	volume={299},
	date={1987},
	pages={431--472}
}
\bib{Soly}{article}{
	title={The boundary distortion and extremal problems in certain classes of univalent functions},
	author={A.Yu. Solynin},
	journal={J. Math. Sci.},
	volume={79},
	date={1996},
	pages={1341–-1358 }
}
\bib{Sol}{article}{
	title={Functional inequalities via polarization},
	author={A.Yu. Solynin},
	journal={St. Petersburg Math. J.},
	volume={8},
	date={1997},
	pages={1015--1038}
}
\bib{Ts}{article}{
	title={A theorem on the majoration of harmonic measure and its applications},
	author={M. Tsuji},
	journal={Tohoku Math. J. (2)},
	volume={3},
	date={1951},
	pages={13--23}
}
\bib{Tsu}{book}{
title={Potential Theory in Modern Function Theory},
author={M. Tsuji},
date={1959},
publisher={Maruzen},
address={Tokyo}
}

\end{biblist}
\end{bibdiv}
\end{document}